# FUNDAMENTAL EUCLIDEAN PATHWISE MINIMIZING EIGENPROPERTIES

DMITRY ZHIGALOV

ABSTRACT. The paper discovers the family of identically-derived Euclidean one-parameter even-dimensional differential linear operators with unique eigenproperties, which prove to be inherently related to the emergent characterizations of fundamental building blocks of embedded minimal surfaces and the Nitsche conjecture proof.

## 1. NTRODUCTION

A while ago a breakthrough in understanding of the minimal surfaces was made through a series of works by T.H. Colding and W.P. Minicozzi, the results summarized in [1], which prove that any embedded minimal surface with finite genus can be built out of planes, helicoids and catenoids – the fundamental building blocks. Moreover, any complete embedded minimal surface with finite topology and one end is either a plane or it has infinite total curvature and its end is asymptotic to a helicoid, while in case of two or more ends, it has finite total curvature and each end thereof is asymptotic to either a plane or a catenoid [2].

What is so very special about the plane, the helicoid and the catenoid? What makes them figure so prominently among the other embedded minimal surfaces with finite topology? What do the ruled minimal surfaces and the catenoid share in common? The fundamental domain of the first Scherk surface is also foliated of catenaries, ones of equal strength, but, unlike catenoid, is not a fundamental building block of embedded minimal surfaces and has planar ends, while the catenoid has catenoidal ends. Why is this? What makes the catenoid and its generatrix, the classic catenary, so special?

Neither geometry nor topology of surfaces can give rigorous answers to the questions arisen above, as the ones proved to be hidden in the fundamental structure of the Euclidean space with, what is remarkable, the classic approach to elicit them. The eigenproperties of differential operators discovered are mutually identical across the entire family. Hence, the least involved case, which corresponds to the Euclidean 3-space, is the wisest to evaluate. Although all the computations presented in the paper are purely analytical, the extensive use of the symbolic computation software was a necessity, as some of them are so bulky that are practically impossible to obtain ‘by hand’ for any dimension of the operator.

## 2. ONE-PARAMETER LINEAR DIFFERENTIAL OPERATOR

Consider the two 3-tuples $\mathbf{u}=(u_1,u_2,u_3)$ and $\mathbf{v}=(v_1,v_2,v_3)$, associated with an orthonormal basis in $\mathbb{R}^3$, $u_i=u_i(t)$, $v_i=v_i(t)$, $u_i\neq 0$, $v_i\neq 0$, $u_i,v_i,t\in\mathbb{R}$, $(i=1,2,3)$. For $t$ some given interval will always be implied with the functions $u_i=u_i(t)$ and $v_i=v_i(t)$ to be continuously differentiable over the entire domain as many times as required. Now suppose $*$ is the operator and $\mathbf{x}$ is its output but such that

$$\mathbf{x}=\mathbf{u}*\mathbf{v}=\begin{bmatrix}u_1v_2\\u_3v_2\\u_3v_1\\u_2v_1\\u_2v_3\\u_1v_3\end{bmatrix} \tag{2.1}$$

Here in (2.1) there is no difference in the particular way of the row-wise arrangement of the diagonal products $u_i v_j$ $\left(i, j = 1,2,3, i \neq j\right)$, as the rearrangement thereof has no effect on the follow-up conclusions. Associated with an orthonormal basis in $\mathbb{R}^6$, unlike the introduced 3-tuples, the 6-tuple $\mathbf{x}$ is *not a vector*, as it does not satisfy the transformation law for its components according to the definition. Its differentiation with respect to $u_i$, $v_i$ $\left(i = 1,2,3\right)$ yields

$$\mathbf{x}\nabla^{\mathrm{T}} = \begin{bmatrix} v_2 & 0 & 0 & 0 & u_1 & 0 \\ 0 & 0 & v_2 & 0 & u_3 & 0 \\ 0 & 0 & v_1 & u_3 & 0 & 0 \\ 0 & v_1 & 0 & u_2 & 0 & 0 \\ 0 & v_3 & 0 & 0 & 0 & u_2 \\ v_3 & 0 & 0 & 0 & 0 & u_1 \end{bmatrix}, \ rank(\mathbf{x}\nabla^{\mathrm{T}}) = 5,$$

Jacobian $\mathbf{x}\nabla^{\mathrm{T}}$ has the property that *any* five of its rows are linearly independent, so elements of any row of the above matrix can be substituted with the unit basis vectors $\boldsymbol{i}_m$ $\left(m = 1,\ldots,6\right)$ to subsequently compute an external product. For the lowest row this gives

$$\mathfrak{D} = \begin{bmatrix} 0 & 0 & v_2 & 0 & u_3 & 0 \\ 0 & 0 & v_1 & u_3 & 0 & 0 \\ 0 & v_1 & 0 & u_2 & 0 & 0 \\ 0 & v_3 & 0 & 0 & 0 & u_2 \\ v_3 & 0 & 0 & 0 & 0 & u_1 \\ \boldsymbol{i}_1 & \boldsymbol{i}_2 & \boldsymbol{i}_3 & \boldsymbol{i}_4 & \boldsymbol{i}_5 & \boldsymbol{i}_6 \end{bmatrix}, \det(\mathfrak{D}) = \begin{bmatrix} u_1 u_2 u_3 v_1 v_3 \\ u_2^2 u_3 v_1 v_3 \\ u_2 u_3^2 v_1 v_3 \\ -u_2 u_3 v_1^2 v_3 \\ -u_2 u_3 v_1 v_2 v_3 \\ -u_2 u_3 v_1 v_3^2 \end{bmatrix}, \tag{2.2}$$

where $\det(\mathfrak{D})$ does not meet the definition of vector as well. Replacement of elements of any other row of $\mathbf{x}\nabla^{\mathrm{T}}$ with $\boldsymbol{i}_m$ $\left(m = 1,\ldots,6\right)$ will result in $\det(\mathfrak{D})$ differing from that of (2.2) only by a factor, which has no effect on the final conclusions (the proof thereof will not be presented in this paper).

Generally, (2.1) can be operation on $k$ $n$-tuples with its diagonal products comprised of $k$ factors, whereas $rank\left(\mathbf{x}\nabla^{\mathrm{T}}\right) = k\left(n-1\right)+1$ if $k < n$, and for $k = n$ the rank is the same as for $k = n-1$. The resulting Euclidean $kn$-space has $k-1$ more dimensions than the number of linearly independent rows of $\mathbf{x}\nabla^{\mathrm{T}}$. Hence, for $k > 2$ external product, as that of (2.2), can not be defined. Similarly, for $n = 2$ there are two 4-tuples.

Differentiation of $\det(\mathfrak{D})$ yields $\det(\mathfrak{D})\nabla^{\mathrm{T}} = u_2 u_3 v_1 v_3 \mathbf{H}$, where

$$\mathbf{H} = \begin{bmatrix} 1 & \frac{u_1}{u_2} & \frac{u_1}{u_3} & \frac{u_1}{v_1} & 0 & \frac{u_1}{v_3} \\ 0 & 2 & \frac{u_2}{u_3} & \frac{u_2}{v_1} & 0 & \frac{u_2}{v_3} \\ 0 & \frac{u_3}{u_2} & 2 & \frac{u_3}{v_1} & 0 & \frac{u_3}{v_3} \\ 0 & -\frac{v_1}{u_2} & -\frac{v_1}{u_3} & -2 & 0 & -\frac{v_1}{v_3} \\ 0 & -\frac{v_2}{u_2} & -\frac{v_2}{u_3} & -\frac{v_2}{v_1} & -1 & -\frac{v_2}{v_3} \\ 0 & -\frac{v_3}{u_2} & -\frac{v_3}{u_3} & -\frac{v_3}{v_1} & 0 & -2 \end{bmatrix} \tag{2.3}$$

Further, the linear map $T:\mathfrak{w}\to\mathfrak{w}'$ is introduced, where $\mathfrak{w}=\mathfrak{w}(t)$ is continuously differentiable as many times as required vector-valued function, and there are such non-zero $u_i$, $v_i$ $(i=1,2,3)$, defined over its entire domain, that

$$\mathfrak{w}'_t=\mathbf{H}\mathfrak{w} \tag{2.4}$$

Let $\mathfrak{w}$ be defined at some fixed $t_0$, then $\mathfrak{w}_0=\mathfrak{w}(t_0)$,

$$\mathfrak{w}_0=\begin{bmatrix}\mathfrak{w}_{0_1}\\ \mathfrak{w}_{0_2}\\ \mathfrak{w}_{0_3}\\ \mathfrak{w}_{0_4}\\ \mathfrak{w}_{0_5}\\ \mathfrak{w}_{0_6}\end{bmatrix} \text{ and } \mathbf{H}\mathfrak{w}_0=\begin{bmatrix}\mathfrak{w}_{0_1}+\frac{u_1\mathfrak{w}_{0_2}}{u_2}+\frac{u_1\mathfrak{w}_{0_3}}{u_3}+\frac{u_1\mathfrak{w}_{0_4}}{v_1}+\frac{u_1\mathfrak{w}_{0_6}}{v_3}\\ 2\mathfrak{w}_{0_2}+\frac{u_2\mathfrak{w}_{0_3}}{u_3}+\frac{u_2\mathfrak{w}_{0_4}}{v_1}+\frac{u_2\mathfrak{w}_{0_6}}{v_3}\\ \frac{u_3\mathfrak{w}_{0_2}}{u_2}+2\mathfrak{w}_{0_3}+\frac{u_3\mathfrak{w}_{0_4}}{v_1}+\frac{u_3\mathfrak{w}_{0_6}}{v_3}\\ -\frac{v_1\mathfrak{w}_{0_2}}{u_2}-\frac{v_1\mathfrak{w}_{0_3}}{u_3}-2\mathfrak{w}_{0_4}-\frac{v_1\mathfrak{w}_{0_6}}{v_3}\\ -\frac{v_2\mathfrak{w}_{0_2}}{u_2}-\frac{v_2\mathfrak{w}_{0_3}}{u_3}-\frac{v_2\mathfrak{w}_{0_4}}{v_1}-\mathfrak{w}_{0_5}-\frac{v_2\mathfrak{w}_{0_6}}{v_3}\\ -\frac{v_3\mathfrak{w}_{0_2}}{u_2}-\frac{v_3\mathfrak{w}_{0_3}}{u_3}-\frac{v_3\mathfrak{w}_{0_4}}{v_1}-2\mathfrak{w}_{0_6}\end{bmatrix} \tag{2.5}$$

Since $rank\left(\left(\mathbf{H}\mathfrak{w}_0\right)\nabla^{\mathrm{T}}\right)=5$ (differentiation is with respect to $u_i$, $v_i$ $(i=1,2,3)$), the function $\mathfrak{w}=\mathfrak{w}(t)$ can not be arbitrary. However, for the vector of the direction cosines of $\left(\mathfrak{w}'_t\right)_0$, $\boldsymbol{\mu}(u_i,v_i)=\frac{\mathbf{H}\mathfrak{w}_0}{|\mathbf{H}\mathfrak{w}_0|}$ $(i=1,2,3)$, $rank\left(\boldsymbol{\mu}(u_i,v_i)\nabla^{\mathrm{T}}\right)=5$. Considering the direction cosines $\frac{\mathfrak{w}_j}{\mathfrak{w}}$ $(j=1,\ldots,6)$ of $\mathfrak{w}=\mathfrak{w}(t)$ always satisfy $\sum_{j=1}^{6}\left(\frac{\mathfrak{w}_j}{\mathfrak{w}}\right)^2=1$, $\mathfrak{w}=\mathfrak{w}(t)$ *can have arbitrary hodograph* $\mathcal{L}$.

## 3. Geometric properties of linear operator $\mathbf{H}$

Differentiation of (2.4) with respect to $t$ yields

$$\mathfrak{w}''_t=\left(\mathbf{H}\mathfrak{w}\right)'_t=\mathbf{H}'_t\mathfrak{w}+\mathbf{H}\mathfrak{w}'_t=\left(\mathbf{H}'_t+\mathbf{H}^2\right)\mathfrak{w}, \tag{3.1}$$

where

$$\mathbf{H}^2=\begin{bmatrix}1 & \frac{2u_1}{u_2} & \frac{2u_1}{u_3} & 0 & 0 & 0\\ 0 & 3 & \frac{2u_2}{u_3} & 0 & 0 & 0\\ 0 & \frac{2u_3}{u_2} & 3 & 0 & 0 & 0\\ 0 & 0 & 0 & 3 & 0 & \frac{2v_1}{v_3}\\ 0 & 0 & 0 & \frac{2v_2}{v_1} & 1 & \frac{2v_2}{v_3}\\ 0 & 0 & 0 & \frac{2v_3}{v_1} & 0 & 3\end{bmatrix}, \tag{3.2}$$

and

$$\mathbf{H}_t' = \begin{bmatrix} 0 & \frac{u_1'}{u_2} - \frac{u_1 u_2'}{u_2^{\,2}} & \frac{u_1'}{u_3} - \frac{u_1 u_3'}{u_3^{\,2}} & \frac{u_1'}{v_1} - \frac{u_1 v_1'}{v_1^{\,2}} & 0 & \frac{u_1'}{v_3} - \frac{u_1 v_3'}{v_3^{\,2}} \\ 0 & 0 & \frac{u_2'}{u_3} - \frac{u_2 u_3'}{u_3^{\,2}} & \frac{u_2'}{v_1} - \frac{u_2 v_1'}{v_1^{\,2}} & 0 & \frac{u_2'}{v_3} - \frac{u_2 v_3'}{v_3^{\,2}} \\ 0 & \frac{u_3'}{u_2} - \frac{u_3 u_2'}{u_2^{\,2}} & 0 & \frac{u_3'}{v_1} - \frac{u_3 v_1'}{v_1^{\,2}} & 0 & \frac{u_3'}{v_3} - \frac{u_3 v_3'}{v_3^{\,2}} \\ 0 & -\frac{v_1'}{u_2} + \frac{v_1 u_2'}{u_2^{\,2}} & -\frac{v_1'}{u_3} + \frac{v_1 u_3'}{u_3^{\,2}} & 0 & 0 & -\frac{v_1'}{v_3} + \frac{v_1 v_3'}{v_3^{\,2}} \\ 0 & \frac{v_2'}{u_2} + \frac{v_2 u_2'}{u_2^{\,2}} & -\frac{v_2'}{u_3} + \frac{v_2 u_3'}{u_3^{\,2}} & -\frac{v_2'}{v_1} + \frac{v_2 v_1'}{v_1^{\,2}} & 0 & -\frac{v_2'}{v_3} + \frac{v_2 v_3'}{v_3^{\,2}} \\ 0 & -\frac{v_3'}{u_2} + \frac{v_3 u_2'}{u_2^{\,2}} & -\frac{v_3'}{u_3} + \frac{v_3 u_3'}{u_3^{\,2}} & -\frac{v_3'}{v_1} + \frac{v_3 v_1'}{v_1^{\,2}} & 0 & 0 \end{bmatrix} \tag{3.3}$$

The next step is to find the *common eigenvectors* $\boldsymbol{\omega}_i$ of $\mathbf{H}^2$ and $\mathbf{H}_t'$, satisfying (2.4). Firstly, the following generalized eigenvalue problem needs to be solved

$$\mathbf{H}^2 \boldsymbol{\omega}_i = \lambda_i \mathbf{H}_t' \boldsymbol{\omega}_i \quad (i = 1, \ldots, q), \tag{3.4}$$

Since $rank\left(\mathbf{H}^2\right) = 6$ and $rank\left(\mathbf{H}_t'\right) = 2$, $q = 2$. The eigenvalues $\lambda_i$ are presented by vector $\boldsymbol{\lambda}$,

$$\boldsymbol{\lambda} = \begin{bmatrix} \dfrac{5 u_2 u_3 v_1 v_3}{u_3 v_1 v_3 u_2' + u_3' u_2 v_1 v_3 - v_1' u_2 u_3 v_3 - v_3' u_2 u_3 v_1} \\ -\dfrac{5 u_2 u_3 v_1 v_3}{u_3 v_1 v_3 u_2' + u_3' u_2 v_1 v_3 - v_1' u_2 u_3 v_3 - v_3' u_2 u_3 v_1} \end{bmatrix},$$

and as the eigenvalues are distinct, there are $q$ generalized eigenvectors, satisfying (3.4).

Solving (3.4) yields

$$\boldsymbol{\omega}_1 = \begin{bmatrix} u_1 \\ u_2 \\ u_3 \\ -v_1 \\ -v_2 \\ -v_3 \end{bmatrix},$$

$$\boldsymbol{\omega}_2 = \begin{bmatrix} \frac{5u_1'u_3'}{4u_3} + \frac{5u_2'u_1'}{4u_2} - \frac{5u_1'v_3'}{4v_3} - \frac{5u_1'v_1'}{4v_1} - \frac{7u_1 u_3'^{\,2}}{8u_3^2} - \frac{u_1 u_2' u_3'}{2u_2 u_3} - \frac{7u_1 u_2'^{\,2}}{8u_2^2} + \frac{u_1 u_3' v_3'}{2u_3 v_3} + \frac{u_1 u_2' v_3'}{2u_2 v_3} + \frac{3u_1 v_3'^{\,2}}{8v_3^2} + \frac{u_1 v_1' u_3'}{2u_3 v_1} + \frac{u_1 u_2' v_1'}{2u_2 v_1} - \frac{u_1 v_1' v_3'}{2v_1 v_3} + \frac{3u_1 v_1'^{\,2}}{8v_1^2} \\ -\frac{7u_2 u_3'^{\,2}}{8u_3^{\,2}} + \frac{3u_2' u_3'}{4u_3} + \frac{3u_2'^{\,2}}{8u_2} + \frac{u_2 v_3' u_3'}{2u_3 v_3} - \frac{3u_2' v_3'}{4v_3} + \frac{3u_2 v_3'^{\,2}}{8v_3^{\,2}} + \frac{u_2 v_1' u_3'}{2u_3 v_1} - \frac{3u_2' v_1'}{4v_1} - \frac{u_2 v_3' v_1'}{2v_1 v_3} + \frac{3u_2 v_1'^{\,2}}{8v_1^{\,2}} \\ \frac{3u_3'^{\,2}}{8u_3} + \frac{3u_2' u_3'}{4u_2} - \frac{7u_3 u_2'^{\,2}}{8u_2^{\,2}} - \frac{3v_3' u_3'}{4v_3} + \frac{u_3 u_2' v_3'}{2u_2 v_3} + \frac{3u_3 v_3'^{\,2}}{8v_3^{\,2}} - \frac{3v_1' u_3'}{4v_1} + \frac{u_3 u_2' v_1'}{2u_2 v_1} - \frac{u_3 v_3' v_1'}{2v_1 v_3} + \frac{3u_3 v_1'^{\,2}}{8v_1^{\,2}} \\ \frac{3v_1 u_2'^{\,2}}{8u_2^{\,2}} - \frac{3u_2' v_1'}{4u_2} - \frac{v_1 u_2' u_3'}{2u_2 u_3} + \frac{v_1 u_2' v_3'}{2u_2 v_3} + \frac{3v_1'^{\,2}}{8v_1} - \frac{3v_1' u_3'}{4u_3} + \frac{3v_3' v_1'}{4v_3} + \frac{3v_1 u_3'^{\,2}}{8u_3^{\,2}} + \frac{v_1 v_3' u_3'}{2u_3 v_3} - \frac{7v_1 v_3'^{\,2}}{8v_3^{\,2}} \\ \frac{3v_2 u_2'^{\,2}}{8u_2^{\,2}} + \frac{v_2 u_2' v_1'}{2u_2 v_1} - \frac{5v_2' u_2'}{4u_2} - \frac{v_2 u_2' u_3'}{2u_2 u_3} + \frac{v_2 u_2' v_3'}{2u_2 v_3} - \frac{7v_2 v_1'^{\,2}}{8v_1^{\,2}} + \frac{5v_2' v_1'}{4v_1} + \frac{v_2 v_1' u_3'}{2u_3 v_1} - \frac{v_2 v_1' v_3'}{2v_1 v_3} - \frac{5v_2' u_3'}{4u_3} + \frac{3v_2 u_3'^{\,2}}{8u_3^{\,2}} + \frac{5v_2' v_3'}{4v_3} + \frac{v_2 v_3' u_3'}{2u_3 v_3} - \frac{7v_2 v_3'^{\,2}}{8v_3^{\,2}} \\ \frac{3v_3 u_2'^{\,2}}{8u_2^{\,2}} + \frac{v_3 u_2' v_1'}{2u_2 v_1} - \frac{3u_2' v_3'}{4u_2} - \frac{v_3 u_2' u_3'}{2u_2 u_3} - \frac{7v_3 v_1'^{\,2}}{8v_1^{\,2}} + \frac{3v_3' v_1'}{4v_1} + \frac{v_3 v_1' u_3'}{2u_3 v_1} + \frac{3v_3'^{\,2}}{8v_3} - \frac{3v_3' u_3'}{4u_3} + \frac{3v_3 u_3'^{\,2}}{8u_3^{\,2}} \end{bmatrix}$$

According to (2.2), $\det(\mathfrak{D}) = u_2 u_3 v_1 v_3 \boldsymbol{\omega}_1$. As $\boldsymbol{\omega}_1$ satisfies (2.4),

$$(\boldsymbol{\omega}_1)'_t = \mathbf{H}\boldsymbol{\omega}_1 = \begin{bmatrix} u_1 \\ u_2 \\ u_3 \\ v_1 \\ v_2 \\ v_3 \end{bmatrix} \tag{3.5}$$

Solving (3.5) yields

$$\boldsymbol{\omega}_1 = \begin{bmatrix} c_1 e^t \\ c_2 e^t \\ c_3 e^t \\ -k_1 e^{-t} \\ -k_2 e^{-t} \\ -k_3 e^{-t} \end{bmatrix}, \tag{3.6}$$

where $c_i = const$, $k_i = const$ $(i = 1, 2, 3)$, or $\boldsymbol{\omega}_1 = \mathfrak{c} e^t + \mathfrak{k} e^{-t}$, where

$$\mathfrak{c} = \begin{bmatrix} c_1 \\ c_2 \\ c_3 \\ 0 \\ 0 \\ 0 \end{bmatrix}, \ \mathfrak{k} = \begin{bmatrix} 0 \\ 0 \\ 0 \\ -k_1 \\ -k_2 \\ -k_3 \end{bmatrix}$$

Considering $\mathfrak{c} \cdot \mathfrak{k} = 0$ and $(\boldsymbol{\omega}_1)'_t = \mathfrak{c} e^t - \mathfrak{k} e^{-t}$, the differential $dA$ of the area swept by vector $\boldsymbol{\omega}_1$ is

$$dA = \frac{1}{2} \begin{vmatrix} ce^t & -ke^{-t} \\ ce^t & ke^{-t} \end{vmatrix} dt = ck\,dt, \tag{3.7}$$

where $c = |\mathfrak{c}|$, $k = |\mathfrak{k}|$, while integration of (3.7) returns

$$\frac{A}{ck} = t + z, \ z = const \tag{3.8}$$

Representing $c$ and $k$ as

$$c = ce^{\frac{a}{2}}, \ k = ce^{-\frac{a}{2}}, \tag{3.9}$$

where $a, c \in \mathbb{R}$, and considering $u = 2A$ as well as, without loss of generality, $z = 0$, expressions (3.8) and (3.9) together yield $\frac{u}{c^2} = 2t$. However, $\omega_1{}^2 = c_1^2 e^{2t} + c_2^2 e^{2t} + c_3^2 e^{2t} + k_1^2 e^{-2t} + k_2^2 e^{-2t} + k_3^2 e^{-2t}$, hence

$$\frac{\omega_1{}^2}{2} = c^2 \cosh\left(\frac{u}{c^2} + a\right), \ u \in [-\mathcal{U}, \mathcal{U}] \tag{3.10}$$

The right part of (3.10), as can be seen, represents the catenary.

Since $\mathbf{H}^2 \boldsymbol{\omega}_1 = 5\boldsymbol{\omega}_1$, it is readily seen that $\boldsymbol{\omega}_1$ is a common eigenvector of $\mathbf{H}^2$ and $\mathbf{H}'_t$, but $\boldsymbol{\omega}_2$ is generally not. The differential equations

$$\mathbf{H}^2 \boldsymbol{\omega}_2 = \lambda_{\mathbf{H}^2_i} \boldsymbol{\omega}_2 \quad (i = 1, \dots, 6), \tag{3.11}$$

need to be solved to find (if one exists) such $\boldsymbol{\omega}_2(t)$, that is, similarly to $\boldsymbol{\omega}_1(t)$, a common eigenvector of $\mathbf{H}'_t$ and $\mathbf{H}^2$ over entire domain.

Eigendecomposition of $\mathbf{H}^2$ yields

$$\boldsymbol{\lambda}_{\mathbf{H}^2} = \begin{bmatrix} 5 \\ 5 \\ 1 \\ 1 \\ 1 \\ 1 \end{bmatrix}, \mathbf{V}_{\mathbf{H}^2} = \begin{bmatrix} \frac{u_1}{u_3} & 0 & 0 & 0 & 0 & 1 \\ \frac{u_2}{u_3} & 0 & 0 & 0 & -\frac{u_2}{u_3} & 0 \\ 1 & 0 & 0 & 0 & 1 & 0 \\ 0 & \frac{v_1}{v_3} & -\frac{v_1}{v_3} & 0 & 0 & 0 \\ 0 & \frac{v_2}{v_3} & 0 & 1 & 0 & 0 \\ 0 & 1 & 1 & 0 & 0 & 0 \end{bmatrix} \tag{3.12}$$

Hence, there are only two distinct eigenvalues with multiplicities 2 and 4 and (3.11) is the set of only two equations.

Now, considering the following system of ordinary differential equations in the matrix form

$$\begin{bmatrix} \mathbf{H}^2\boldsymbol{\omega}_2 \\ \mathbf{H}\boldsymbol{\omega}_2 \end{bmatrix} = \begin{bmatrix} \lambda_{\mathbf{H}^2 i}\boldsymbol{\omega}_2 \\ (\boldsymbol{\omega}_2)'_t \end{bmatrix}, \tag{3.13}$$

it is readily seen that (3.13) is ovedetermined.

Solving (3.13) for $\lambda_{\mathbf{H}^2} = 5$ yields the following general solution

$$\boldsymbol{\omega}_2 = \begin{bmatrix} \mathbb{C}_1 e^{5t} \\ \mathbb{C}_2 e^{5t} \\ \mathbb{C}_3 e^{5t} \\ \mathbb{K}_1 e^{-5t} \\ \mathbb{K}_2 e^{-5t} \\ \mathbb{K}_3 e^{-5t} \end{bmatrix} = \begin{bmatrix} u_1 \\ u_2 \\ u_3 \\ v_1 \\ v_2 \\ v_3 \end{bmatrix}, \tag{3.14}$$

where $\mathbb{C}_i = const$ and $\mathbb{K}_i = const$ $(i = 1, 2, 3)$. For the eigenvalue $\lambda_{\mathbf{H}^2} = 1$ system (3.13) has no solutions.

Let $\boldsymbol{\omega}_i = \nu\boldsymbol{\omega}_i$ $(i = 1, 2)$, where $\nu = \nu(t)$ is an arbitrary real function continuously differentiable over entire domain, and $\boldsymbol{\omega}_i$ satisfies (2.4), then $(\boldsymbol{\omega}_i)'_t = \mathbf{H}\boldsymbol{\omega}_i = \nu\mathbf{H}\boldsymbol{\omega}_i = \nu(\boldsymbol{\omega}_i)'_t$, which means that $\nu(t) = const$. Hence, $\boldsymbol{\omega}_i$ has the same form as (3.6) and (3.14), so the latter two *are the general solutions*.

It is easily seen, that $\frac{\boldsymbol{\omega}_2^2}{2}$ is equal to the right part of (3.10). The solution (3.14) is, in fact, *the conjugate of the solution (3.6)*.

Thus, *there are only two vector-valued functions* $\boldsymbol{\omega}_1$ *and* $\boldsymbol{\omega}_2$*, satisfying (2.4), such that each one is a common eigenvector of* $\mathbf{H}^2$ *and* $\mathbf{H}'_t$ *over entire domain, and 1 radian of area swept by each of them as a function of* $t$ *represents the catenary as shown in (3.10).*

Suppose $t = t(u)$ is twice continuously differentiable, where $u$, as was stated above, is the double sweep area of $\mathfrak{w}$. Then $\mathfrak{w}'_u = \mathfrak{w}'_t t'_u$, and from (2.4) it follows

$$\mathfrak{w}'_u = t'_u \mathbf{H}\mathfrak{w} \tag{3.15}$$

Then, considering (3.15),

$$\mathfrak{w}''_u = (t'_u \mathbf{H}\mathfrak{w})'_u = (t'_u \mathbf{H})'_u \mathfrak{w} + (t'_u \mathbf{H})\mathfrak{w}'_u = (t''_u \mathbf{H} + t'_u \mathbf{H}'_t)\mathfrak{w} + (t'_u \mathbf{H})t'_u \mathbf{H}\mathfrak{w} = t''_u \mathfrak{w}'_t + {t'_u}^2 \left(\mathbf{H}'_t\mathfrak{w} + \mathbf{H}^2\mathfrak{w}\right) \tag{3.16}$$

The functions $\mathfrak{w}''_u$ and $\mathfrak{w}(t(u))$ are collinear over entire domain, as their argument is the linear function of the sweep area of $\mathfrak{w}$. Then, as follows from the findings, with the only exceptions to be $\mathfrak{w} = \omega_i$, there are *non-vanishing all at once* orthogonal projections of $t''_u \mathfrak{w}'_t$, $t'^2_u \mathbf{H}'_t \mathfrak{w}$ and $t'^2_u \mathbf{H}^2 \mathfrak{w}$ onto $\mathfrak{w}$ over at least some subset of the domain of $\mathfrak{w}$, but their sum vanishes everywhere. This is what can be called the *measure dissipation* of that subset.

Once the function $t = t(u)$ is linear, the term $t''_u \mathfrak{w}'_t$ vanishes, and, as can be proven, the two remaining terms of (3.16) align with $\mathfrak{w}$.

Let $s(u) = \frac{\mathfrak{w}(t(u))^2}{2}$, then, once $d\mathfrak{w}$ and $d\mathfrak{w}_\tau$ are the absolute values of the radial and the transverse components of $d\mathfrak{w}$ with respect to $\mathfrak{w}$, denoting $d\ell$ the arc length differential of the hodograph $\mathcal{L}$, the arc length differential $dl$ of $s = s(u)$ has the form

$$dl = \sqrt{(ds)^2 + (du)^2} = \sqrt{\left[d\left(\frac{\mathfrak{w}^2}{2}\right)\right]^2 + (\mathfrak{w} d\mathfrak{w}_\tau)^2} = \mathfrak{w}\sqrt{(d\mathfrak{w})^2 + (d\mathfrak{w}_\tau)^2} = \mathfrak{w} d\ell ,$$

then

$$\frac{ds}{dl} = d\left(\frac{\mathfrak{w}^2}{2}\right)\frac{1}{dl} = \frac{\mathfrak{w} d\mathfrak{w}}{\mathfrak{w} d\ell} = \frac{d\mathfrak{w}}{d\ell} \tag{3.17}$$

As $2\pi s = \pi \mathfrak{w}^2$, the function $\mathfrak{w}$ represents the radius of the solid of revolution defined around $l$ – axis. The total volume thereof is equal to the total area of the respective surface of revolution generated by the function $s(u)$. Hence, $\omega_1$ and $\omega_2$ are the radii of the solids of revolution with the total volume equal to the surface area of the respective catenoids.

Given $u = u(l)$ and $l = l(\ell)$ are continuously differentiable, from (3.17) it is seen that $\mathfrak{w}(t(u(l(\ell))))$ and $s(u(l))$ share identical behaviors, which means that the measure dissipation, if present for $\mathfrak{w}(t(u(l(\ell))))$, is fully inherited by $s(u(l))$. That is *why exactly the double sweep area* $u$ *matters*. Dissipation itself, as it is easy to see, is only the attribute of the sweep area of $\mathfrak{w}$, as the radial component of the arc length differential $d\ell$ is by its very nature varies the norm of $\mathfrak{w}$ without dissipation. Rather, it is the domain $[-\mathcal{U}, \mathcal{U}]$ which defines the integral double sweep area of $\mathfrak{w}$ and may be *dissipated*. And if $\mathfrak{w} = \omega_i$ over entire domain, there is no dissipation at all. Hence, the classic catenary is the *non-dissipative* path.

It can be shown that $\mathfrak{w}'_t$, $\mathbf{H}'_t \mathfrak{w}$ and $\mathbf{H}^2 \mathfrak{w}$ are generally *linearly independent* even for the flat hodographs. Hence, the 2-dimensional geometry fails to fully describe planar curves!

As is known, among all minimal graphs over given annulus the upper slab of catenoid has the greatest conformal modulus [3], which is the essence of the Nitsche conjecture, and it is readily obvious that such dissipation is responsible for a lesser value of the conformal modulus. Putting it simply, in case of catenoid, the value $u$ is completely mapped into the energy minimizing shrinking (in terms of soap film physics), and, as catenoid has axial symmetry, this mapping is radially uniform. Hence, the catenoid is the most ‘slack’ of all such minimal graphs, so it can extend itself higher. The measure dissipation occurs for all other possible annular graphs, and the Nitsche conjecture can be proven by such characterization.

Fig. 1 For a fixed length hodograph $\mathcal{L}$ once $u=0$, function $\mathfrak{w}$, extending to the end of $\mathcal{L}$ (both shown in black), has always the least possible norm against the case of $u\neq 0$ (shown in red), given initial value $\mathfrak{w}_0$.

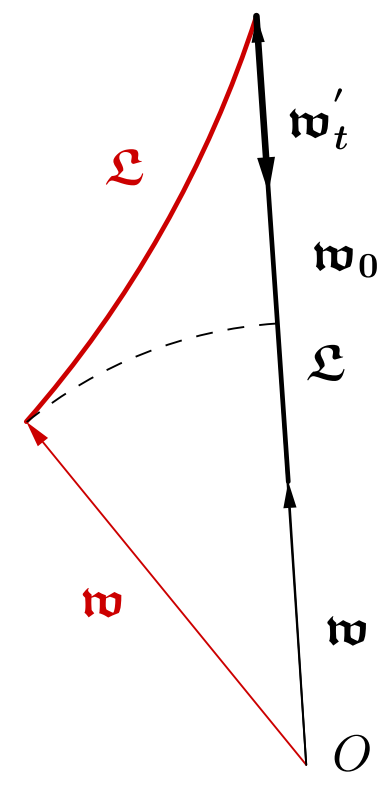

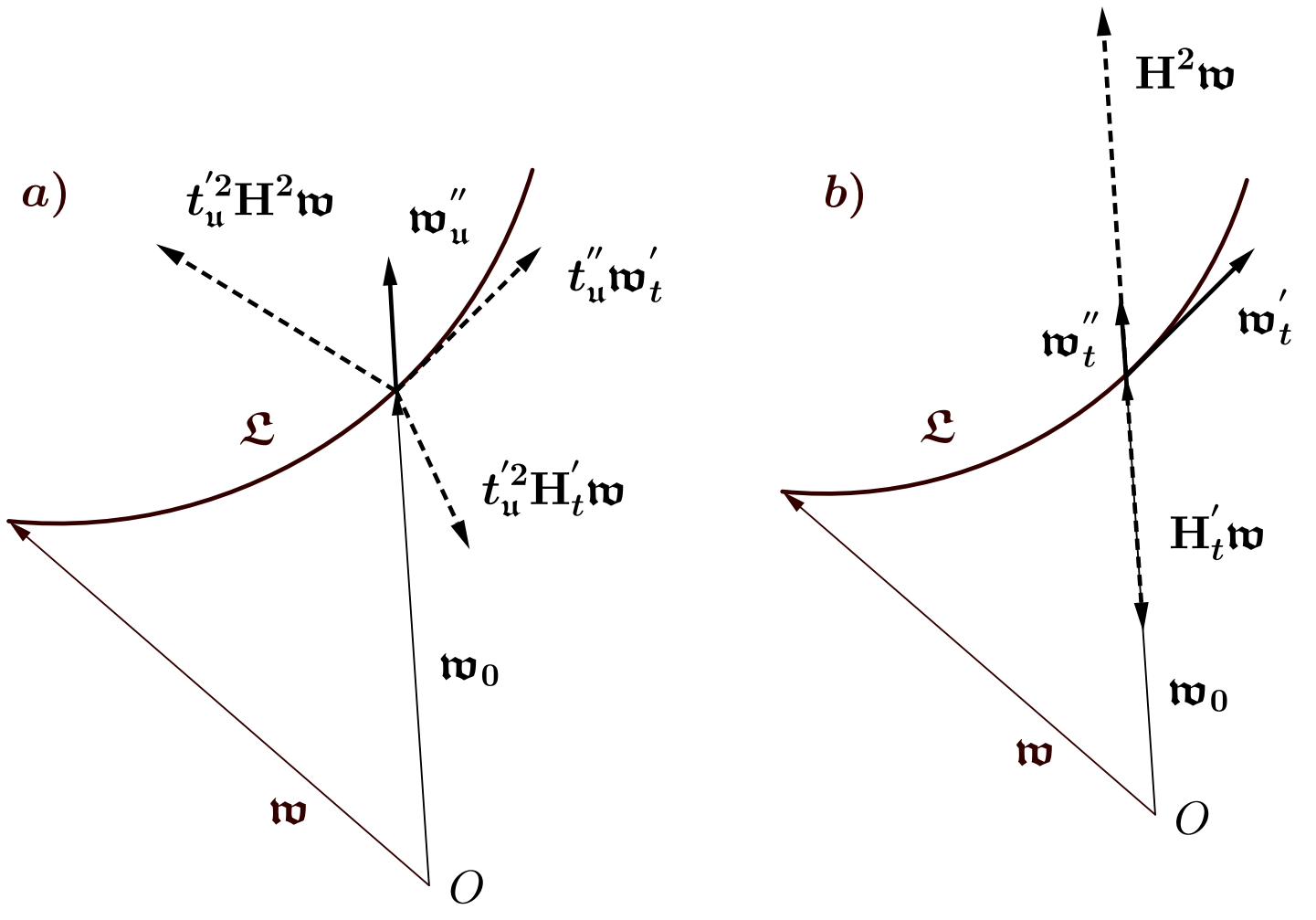

Fig. 2 The case **a** is a general one with $t''_u\mathfrak{w}'_t$, $t'^2_u\mathbf{H}'_t\mathfrak{w}$ and $t'^2_u\mathbf{H}^2\mathfrak{w}$ to be linearly independent, but summing up to the second order derivative directed along the function itself, as the argument thereof is its own double sweep area. The case **b** is the alignment of $\mathbf{H}'_t\mathfrak{w}$ and $\mathbf{H}^2\mathfrak{w}$ with the function, which corresponds to zero measure dissipation. In this case, the initial argument, as derived from solution of (2.4), is the linear function of the double sweep area $u$. Hence, the terms emerge with no multipliers.

Equation (3.10) does not always have a solution and, as is known from physics of soap films, in such cases, the film collapses and forms flat areas around the rings, which had been the boundaries of the catenoid. This case is known as the Goldschmidt solution. However, there is no sweep area now, i.e. $\mathcal{U}=0$, and no dissipation as the result since the norm change of $\mathfrak{w}$ with each $d\ell$ is $d\ell$ itself. Hence, the Goldschmidt solutions are foliations of coplanar concentric circles shrinking to their common center. In this case, the non-dissipative path is the straight line.

Figs. 1-2 describe the possible pathwise behaviors related to the measure dissipation.

Thus, the minimizing paths are associated with, what it appears to be, the *inner structure* of the derivatives of $\mathfrak{w}$. While $\mathfrak{w}'_t$ has no inner structure, the second order derivative $\mathfrak{w}''_t$ is split in, generally, three components, which are derived with the differential split-operator $\mathbf{H}$, and when all the components align with the function itself, the minimizing paths occur.

But are there similar ‘dissipation-free’ cases possible for higher order derivatives of $\mathfrak{w}$? Should they exist, there must be other minimizing paths in addition to the given two. Differentiation of (3.1) yields 4 linear operators with the void set of generalized eigenvectors for two of them, $\left(\mathbf{H}^2\right)'_t$ and $\mathbf{H}''_t$.

Differentiating the third order derivative $p$ times $(p = 1, 2 \ldots)$, the operators $\left(\mathbf{H}^2\right)_t' \mathbf{H}_t^p$ and $\mathbf{H}_t'' \mathbf{H}_t^p$ are generated. Then, let the following generalized eigenvalue problem be introduced

$$\left(\mathbf{H}^2\right)_t' \mathbf{H}_t^p \mathbf{z} = \alpha \mathbf{H}_t'' \mathbf{H}_t^p \mathbf{z} \, ,$$

where $\mathbf{z}$ is a generalized eigenvector and $\alpha$ is the eigenvalue. It can be reduced to

$$\left(\mathbf{H}^2\right)_t' \mathbf{f} = \alpha \mathbf{H}_t'' \mathbf{f} \, , \tag{3.18}$$

where $\mathbf{f} = \mathbf{H}_t^p \mathbf{z}$, which means $\mathbf{f}$ is not arbitrary. However, (3.18) does not have a solution anyway since $\left(\mathbf{H}^2\right)_t'$ and $\mathbf{H}_t''$ do not have generalized eigenvectors. Hence, starting from $\mathfrak{w}_t'''$ there is no common eigenvector for all the sub-operator addends generated by the subsequent differentiation of (2.4) and splitting the derivative.

Thus, the line and the (classic) catenary are the only possible non-dissipative paths in Euclidean geometry. Considering this, the plane, the helicoid and the catenoid can be characterized *as foliations of cylindrical (spherical) geodesics which are pathwise-connected with non-dissipative paths, orthogonal everywhere to the leaves of foliation, in such way, that their radius is varied according to the non-dissipative characterizations, with the geodesics itself shrinking to the lowest possible limit of its length*. For a circle the non-dissipative circumference variation is directly dependent on that of its radius, while for a helix the arc length is the scaled length of the hypotenuse, with the radius to be one of the catheti (another one, defining pitch, is fixed). Hence, once the radius value is varied in non-dissipative way, the helix arc length is varied so as well. The minimal surfaces which are such foliations can be called *perfect* minimal surfaces. They ‘emulate’ the minimizing behavior of the non-dissipative paths.

All minimal surfaces are characterized as having zero mean curvature everywhere, while perfect minimal surfaces are also characterized with zero measure dissipation everywhere. It is remarkable that the ends of complete embedded minimal surfaces with finite topology have only three possible asymptotic behaviors, namely, planar, helicoidal and catenoidal, and the perfect minimal surfaces are the only fundamental building blocks of those surfaces.

## 4. Geometric properties of linear operator $\mathbf{H}$

Eigendecomposition of linear operator $\mathbf{H}$ gives

$$\boldsymbol{\lambda}_{\mathbf{H}}=\begin{bmatrix}\sqrt{5}\\-\sqrt{5}\\1\\1\\-1\\-1\end{bmatrix},\ \mathbf{V}_{\mathbf{H}}=\begin{bmatrix}-\frac{1}{2}\frac{u_1\left(\sqrt{5}-1\right)}{\left(\sqrt{5}-2\right)} & -\frac{1}{2}\frac{u_1\left(-\sqrt{5}-1\right)}{\left(-\sqrt{5}-2\right)} & 0 & 1 & 0 & 0\\ -\frac{1}{2}\frac{u_2\left(\sqrt{5}-1\right)}{\left(\sqrt{5}-2\right)} & -\frac{1}{2}\frac{u_2\left(-\sqrt{5}-1\right)}{\left(-\sqrt{5}-2\right)} & -u_2 & 0 & 0 & 0\\ -\frac{1}{2}u_3\left(3+\sqrt{5}\right) & -\frac{1}{2}u_3\left(3-\sqrt{5}\right) & u_3 & 0 & 0 & 0\\ v_1 & v_1 & 0 & 0 & 0 & -v_1\\ v_2 & v_2 & 0 & 0 & 1 & 0\\ v_3 & v_3 & 0 & 0 & 0 & v_3\end{bmatrix}$$

Singular value decomposition of $\mathbf{H}$ yields vector $\boldsymbol{\sigma}_{\mathbf{H}}$ of the singular values, which, as it will be seen further, is always defined in the form

$$\boldsymbol{\sigma}_{\mathbf{H}}=\begin{bmatrix}\sigma_{\mathbf{H}\max}\\1\\1\\1\\1\\\sigma_{\mathbf{H}\min}\end{bmatrix}, \tag{4.1}$$

where

$$\begin{aligned}\sigma_{\mathbf{H}\max}&=\frac{\sqrt{2}}{2}\left(\frac{1}{u_2^2u_3^2v_1^2v_3^2}\left(b+q\right)\right)^{\frac{1}{2}}\\ \sigma_{\mathbf{H}\min}&=\frac{\sqrt{2}}{2}\left(\frac{1}{u_2^2u_3^2v_1^2v_3^2}\left(b-q\right)\right)^{\frac{1}{2}}\end{aligned}, \tag{4.2}$$

and

$$b=\left(\left(v_1^2+v_3^2\right)u_3^2+v_1^2v_3^2\right)u_2^4+u_3^4v_1^2v_3^2+\left(v_3^2v_1^4+\left(v_3^4+\left(u_1^2+v_2^2\right)v_3^2\right)v_1^2\right)u_3^2+$$
$$+\left(\left(v_1^2+v_3^2\right)u_3^4+\left(v_1^4+\left(u_1^2+v_2^2+14v_3^2\right)v_1^2+v_3^4+\left(u_1^2+v_2^2\right)v_3^2\right)u_3^2+v_3^2v_1^4+\left(v_3^4+\left(u_1^2+v_2^2\right)v_3^2\right)v_1^2\right)u_2^2,$$
$$q=\left(\left(u_1^2+u_2^2+u_3^2+v_1^2+v_2^2+v_3^2\right)\left(\left(\left(v_1^2+v_3^2\right)u_3^2+v_1^2v_3^2\right)u_2^2+u_3^2v_1^2v_3^2\right)\left(\left(\left(v_1^2+v_3^2\right)u_3^2+v_1^2v_3^2\right)u_2^4+\right.\right.$$
$$\left.\left.+\left(\left(v_1^2+v_3^2\right)u_3^4+\left(v_1^4+\left(u_1^2+v_2^2+24v_3^2\right)v_1^2+v_3^2\left(u_1^2+v_2^2+v_3^2\right)\right)u_3^2+v_1^2v_3^2\left(u_1^2+v_1^2+v_2^2+v_3^2\right)\right)u_2^2+u_3^2v_1^2v_3^2\left(u_1^2+u_3^2+v_1^2+v_2^2+v_3^2\right)\right)\right)^{\frac{1}{2}}$$

According to (4.2) and the relation between determinant and singular values, $\left|\det\left(\mathbf{H}\right)\right|=\sigma_{\mathbf{H}\max}\sigma_{\mathbf{H}\min}=5$.

Now the following matrix is to be introduced

$$\mathbf{N}=\begin{bmatrix} u_1 & u_1 & u_1 & -u_1 & u_1 \\ -u_2 & u_2 & -u_2 & -u_2 & -u_2 \\ u_3 & u_3 & u_3 & u_3 & u_3 \\ v_1 & -v_1 & -v_1 & -v_1 & -v_1 \\ -v_2 & -v_2 & v_2 & v_2 & -v_2 \\ -v_3 & -v_3 & v_3 & v_3 & v_3 \end{bmatrix},$$

where it is seen that $rank(\mathrm{N})=5$ for any $u_i$, $v_i$, $u_i \neq 0$, $v_i \neq 0$ $(i=1,2,3)$. Then the linear combination $\mathfrak{e}$ of columns $\mathbf{N}_j$ $(j=1,\dots,6)$ of matrix $\mathbf{N}$ is

$$\mathfrak{e}=\sum_{j=1}^{6}\varpi_j\mathbf{N}_j=\begin{bmatrix} u_1(-\varpi_4+\varpi_5+\varpi_3+\varpi_2+\varpi_1) \\ -u_2(\varpi_1-\varpi_2+\varpi_3+\varpi_4+\varpi_5) \\ u_3(\varpi_1+\varpi_2+\varpi_3+\varpi_4+\varpi_5) \\ v_1(\varpi_1-\varpi_2-\varpi_3-\varpi_4-\varpi_5) \\ -v_2(-\varpi_4+\varpi_5-\varpi_3+\varpi_2+\varpi_1) \\ -v_3(\varpi_1+\varpi_2-\varpi_3-\varpi_4-\varpi_5) \end{bmatrix},$$

where $\varpi_j$ $(j=1,\dots,6)$ are arbitrary non-zero real values.

Denoting $\mathbf{p}=\mathbf{H}\mathfrak{e}$,

$$\mathbf{p}=\begin{bmatrix} u_1(-\varpi_4+\varpi_5+\varpi_3+\varpi_2+\varpi_1) \\ -u_2(\varpi_1-\varpi_2+\varpi_3+\varpi_4+\varpi_5) \\ u_3(\varpi_1+\varpi_2+\varpi_3+\varpi_4+\varpi_5) \\ -v_1(\varpi_1-\varpi_2-\varpi_3-\varpi_4-\varpi_5) \\ v_2(-\varpi_4+\varpi_5-\varpi_3+\varpi_2+\varpi_1) \\ v_3(\varpi_1+\varpi_2-\varpi_3-\varpi_4-\varpi_5) \end{bmatrix}$$

Computation gives $|\mathfrak{e}|=|\mathbf{p}|$. Hence, the circular section of the hyperellipsoid generated by $\mathbf{H}$ in $\mathbb{R}^6$ is the unit hypersphere, defined by matrix $\mathbf{N}$.

The linear combination of the two first columns of $\mathbf{V}_{\mathbf{H}^2}$, defined in (3.12), gives

$$\boldsymbol{\delta}=\gamma_1\mathbf{v}_{\mathbf{H}^2 1}+\gamma_2\mathbf{v}_{\mathbf{H}^2 2}=\begin{bmatrix} \frac{\gamma_1 u_1}{u_3} \\ \frac{\gamma_1 u_2}{u_3} \\ \gamma_1 \\ \frac{\gamma_2 v_1}{v_3} \\ \frac{\gamma_2 v_2}{v_3} \\ \gamma_2 \end{bmatrix} \quad (i=1,\dots,6), \tag{4.3}$$

where $\gamma_1$ and $\gamma_2$ are arbitrary non-zero real values.

Expanding $\mathbf{N}$ with $\boldsymbol{\delta}$ as below,

$$\mathbf{E}=\begin{bmatrix} u_1 & u_1 & u_1 & -u_1 & u_1 & \frac{\gamma_1 u_1}{u_3} \\ -u_2 & u_2 & -u_2 & -u_2 & -u_2 & \frac{\gamma_1 u_2}{u_3} \\ u_3 & u_3 & u_3 & u_3 & u_3 & \gamma_1 \\ v_1 & -v_1 & -v_1 & -v_1 & -v_1 & \frac{\gamma_2 v_1}{v_3} \\ -v_2 & -v_2 & v_2 & v_2 & -v_2 & \frac{\gamma_2 v_2}{v_3} \\ -v_3 & -v_3 & v_3 & v_3 & v_3 & \gamma_2 \end{bmatrix},$$

and

$$\det(\mathbf{E})=-32(\gamma_2 u_3+\gamma_1 v_3)v_2 v_1 u_2 u_1$$

Once the condition $\gamma_2 u_3+\gamma_1 v_3=0$ holds, $\det(\mathbf{E})=0$, which implies

$$\gamma_1=\theta u_3,\ \gamma_2=-\theta v_3, \tag{4.4}$$

where $\theta$ is an arbitrary non-zero real value.

Considering (4.4), expression (4.3) takes form

$$\boldsymbol{\delta}=\theta\begin{bmatrix} u_1 \\ u_2 \\ u_3 \\ -v_1 \\ -v_2 \\ -v_3 \end{bmatrix}$$

Thus, *vector* $\boldsymbol{\delta}$ *is always the intersection of the circular cross section of the hyperellipsoid of* $\mathbf{H}$ *and the eigenspace* $E_5(\mathbf{V}_{\mathbf{H}^2})$ *of the eigenvalue* $\lambda_{\mathbf{H}^2}=5$ *of the matrix* $\mathbf{H}^2$. Considering this, $\sigma_{\mathbf{H}\max}>1, \sigma_{\mathbf{H}\min}<1$. Hence, the vector $\boldsymbol{\sigma}_{\mathbf{H}}$ is always defined in the form (4.1).

## References

[1] T.H. Colding and W.P. Minicozzi II, Shapes of embedded minimal surfaces, Proc. Natl. Acad. Sci. USA 103 (2006), no. 30, 11106–11111.

[2] J. Bernstein and C. Breiner. Conformal structure of minimal surfaces with finite topology. Comm. Math. Helv., 86(2):353–381, 2011. MR2775132, Zbl 1213.53011.

[3] Iwaniec, T., Kovalev, L.V., Onninen, J.: The Nitsche conjecture. J. Am. Math. Soc. 24(2), 345–373 (2011)